\documentclass[preprint,12pt]{elsarticle}

\usepackage{amsmath,amssymb,amsfonts}
\usepackage{graphicx}
\usepackage[notcite,notref]{showkeys}
\usepackage{graphicx,color}

\newtheorem{theorem}{Theorem}[section]
\newtheorem{definition}[theorem]{Definition}
\newtheorem{proposition}[theorem]{Proposition}
\newtheorem{corollary}[theorem]{Corollary}

\newtheorem{lemma}[theorem]{Lemma}

\def\BT{{\mathbb T}}

\newcommand{\ands}{\quad\mbox{and}\quad}

\def\iy{\infty}
\def\ts{\times}
\def\wt{\widetilde}

\addcontentsline{toc}{section}{text_to_be_added}

\journal{Indagationes Math.}

\begin{document}
\begin{frontmatter}
\title{Wiener-Hopf factorization indices of rational matrix functions
with respect to the unit circle in terms of realization}
\author{G.J. Groenewald}
\ead{gilbert.groenewald@nwu.ac.za}
\address{Gilbert J. Groenewald, School~of~Mathematical~and Statistical Sciences, 
North-West~University, 
Research Focus: Pure and Applied Analytics, 
Private~Bag~X6001, 
Potchefstroom~2520, 
South Africa  }
\author{M.A.Kaashoek}
\ead{m.a.kaashoek@vu.nl}
\address{Marinus A. Kaashoek, Afdeling Wiskunde, Faculteit
    der Exacte Wetenschappen, Vrije Universiteit Amsterdam, De Boelelaan
    1111, 1081 HV Amsterdam, The Netherlands}
\author{A.C.M. Ran}
\ead{a.c.m.ran@vu.nl}
\address{Andr\'e C. M. Ran,
Department of Mathematics, Faculty of Sciences, Vrije Universiteit Amsterdam, De Boelelaan
    1111, 1081 HV Amsterdam, The Netherlands and Research Focus: Pure and Applied Analytics, North West University, Potchefstroom, South Africa}

\begin{abstract}
As in the paper \cite{GKR} our  aim is to obtain explicitly  the Wiener-Hopf indices  of  a rational $m\ts m$ matrix function  $R(z)$ that has no poles and no zeros on the unit circle $\BT$ but, in contrast with \cite{GKR}, the function $R(z)$ is not required to be unitary on the unit circle.   On the other hand, using  a Douglas-Shapiro-Shields type of factorization, we  show that $R(z)$  factors as $R(z)=\Xi(z)\Psi(z)$, where $\Xi(z)$ and $\Psi(z)$ are rational $m\ts m$ matrix functions, $\Xi(z)$ is unitary  on the unit circle and $\Psi(z)$ is an invertible outer function.   Furthermore,  the fact that $\Xi(z)$ is unitary  on the unit circle allows us  to factor  as $\Xi(z) =V(z)W^*(z)$ where $V(z)$ and $W(z)$ are rational bi-inner $m\ts m$ matrix functions. The latter    allows us  to  solve the Wiener-Hopf indices problem. To derive explicit formulas  for  the functions $V(z)$ and $W(z)$  requires additional realization properties of the function  $\Xi(z)$ which are given in the last two sections. 
\end{abstract}
\end{frontmatter}
\begin{center} 
\emph{Dedicated to Jaap Korevaar, on the occasion of his\\ 100-th birthday, in friendship and with gratitude.}
\end{center}

\setcounter{equation}{0}
\section{Introduction}\label{sec1:intro}

Throughout this paper $R(z)$ is an $m\ts m$ rational matrix function which has no poles and no zeros   on the unit circle $\BT$.  Such a function admits  a right Wiener-Hopf factorization, (see, for instance, \cite{BGKOT21, BGKROT200, CGOT3}), that is, 
\[
R(z)=W_-(z)D(z)W_+(z),
\]
where the factors $W_-(z)$ and $W_+(z)$ are rational $m\times m$ matrix functions such that
\begin{itemize}
\item[{\rm (i)}] $W_-(z)$ has no poles and no zeros outside the open unit disc including infinity,
\item[{\rm (ii)}] $W_+(z)$ has no poles and no zeros on the closed unit disc,
\end{itemize}
and where the middle term $D(z)$ is a diagonal matrix 
$$
D(z)={\rm diag\, }(z^{-\alpha_1}, \cdots , z^{-\alpha_s}) \oplus I_k \oplus {\rm diag\, }
(z^{\omega_t}, \cdots , z^{\omega_1}),
$$
where $-\alpha_1\leq \cdots \leq -\alpha_s<0$ and $0<\omega_t\leq \cdots\leq \omega_1$ are integers,
and $m=s+k+t$. The numbers $\alpha_j$ and $\omega_j$ are called the \emph{right Wiener-Hopf indices of $R(z)$.} Reversing the roles of $W_-(z)$ and $W_+(z)$ one obtains the definition of a \emph{left Wiener-Hopf factorization} and \emph{left Wiener-Hopf indices.} In what follows we will restrict to the right Wiener-Hopf indices and omit the word ``right''. 

Wiener-Hopf factorization of (rational) matrix functions plays a role in the study of Toeplitz operators, convolution integral operators and singular integral operators, in particular, in the study of the Fredholm properties of such operators. See, e.g., \cite{CGOT3}, \cite{GF}, \cite{GGKOT49}, Chapters XII and XIII, and \cite{GGKOT63}, Chapter XXIV. For more details, see also the introduction to \cite{GKR}, and the references given there.

The main goal of the present paper is to obtain explicitly the Wiener-Hopf indices of  any  $m\ts m$ rational matrix function $R(z)$ which has no poles and no zeros  on the unit circle $\BT$. We call this the  \emph{Wiener-Hopf indices problem}. This problem has been solved in \cite{GKR} for the  case when additionally  $R(z)$ is unitary for any $z\in \BT$.  In the present paper   we shall solve the Wiener-Hopf indices problem for any $R(z)$ not necessarily unitary on $\BT$. The analysis is based on the fact that the function $R(z)$ has a realization of the following type:
\begin{equation}\label{R}
R(z)=R_0+zC(I-zA)^{-1}B+ \gamma (zI-\alpha)^{-1}\beta,
\end{equation}
with the square matrices $A$ and $\alpha$ being stable, i.e., the eigenvalues of $A$ and $\alpha$ are in the open unit disc. Such a realization for $R(z)$, with $A$ and $\alpha$ stable, exists because $R(z)$ has no poles on the unit circle $\mathbb{T}$, see, e.g., \cite{BGKROT178}, in particular Chapter 8.

It is the intention of this paper to solve the Wiener-Hopf indices problem by carrying out all the  steps explicitly, starting from the  realization \eqref{R}. In effect the main aim is finding an algorithm for computing the Wiener-Hopf indices of $R(z)$.

The paper consists of six sections including the present introduction. In the second section, using the realization \eqref{R},  an analogous realization  is presented for the product  $R^*(z)R(z)$, where $R^*(z)=R(1/\bar{z})^*$, and this realization is used  to obtain in realized form an invertible outer factor $\Psi(z)$  such that
\begin{equation}\label{RstarR=PsistarPsi}
R^*(z)R(z)=\Psi^*(z)\Psi(z).
\end{equation}
Recall that a rational matrix function $\Psi(z)$ is called an \emph{invertible outer function} if $\Psi(z)$ has no poles and no zeros on the closed unit disc. Moreover, an invertible outer factor of $R(z)$ is unique up to multiplication on the left by a constant unitary matrix (see, e.g., \cite{FB}, Theorems 5.2.1 and 6.1.1). 
The realization for $\Psi(z)$ is given in Proposition \ref{prop:realpsi}. The matrices in the realization can be constructed explicitly fromt he matrices in the realization of $R(z)$, the construction involves the solutions to two Stein equations and a discrete algebraic Riccati equation.

In the third section we introduce the function 
\begin{equation}\label{definitionXi}
\Xi (z)=R(z)\Psi(z)^{-1}.
\end{equation}
Obviously, this function is an $m\ts m$ rational matrix function.
Notice that $\Xi(z)$ is uniquely determined by $\Psi(z)$, and hence is unique up to multiplication by a constant unitary matrix on the right.
Lemma \ref{lem:Xi} shows that $\Xi(z)$ is unitary for each $z\in \BT$ and   the functions $R(z)$ and $\Xi(z)$  have the same right Wiener-Hopf  indices. 
Thus in order to solve the Wiener-Hopf  indices problem   it suffices to solve this problem for the function  $\Xi(z)$.  Since $\Xi(z)$ is unitary for each $z\in \BT$,  we can apply Theorem 1.1 of  \cite{GKR} to get these indices. To get the required  explicit formulas we need to factorize  $\Xi(z)$ as $\Xi(z)=V(z)W^*(z) $ where  $V(z)$ and $W(z)$ are rational bi-inner $m\ts m$ matrix functions with additional results that  are presented in the final two sections. Recall that a rational matrix function is called \emph{bi-inner} if it has unitary values on the unit circle and no poles on the unit disc.
In Section  4 we present $\Xi(z)$ in realized form (see \eqref{realizationXi}) and a few related results.
The realization is based on the realizations for $R(z)$ and $\Psi(z)$, and involves the solution to an additional Stein equation. In Section 5 we construct explicitly the bi-inner functions $V(z)$ and $W(z)$  appearing in the factorization $\Xi(z)=V(z)W(z)^*$. The explicit construction of realizations for $V(z)$ and $W(z)$ is based on an algorithm described in Section 4.7 in \cite{FB}, and involves the solution to yet one more Stein equation and a Lyapunov equation. Finally, in Section 6 we discuss some additional properties of $V(z)$ and $W(z)$ and their realizations.

\setcounter{equation}{0}
\section {The invertible outer factor of  $R^*(z)R(z)$} 
Throughout $R(z)$ is the rational $m\ts m$ matrix function given by \eqref{R}, that is, 
\begin{equation}\label{R2}
R(z)=R_0+zC(I-zA)^{-1}B+ \gamma (zI-\alpha)^{-1}\beta, 
\end{equation}
with the square matrices $A$ and $\alpha$ being stable, i.e., $A$ and $\alpha$ have all their eigenvalues located in the open unit disc.  Furthermore,   $P_+$ and $P_-$ are  the observability gramians given by
\begin{equation}\label{eq:obsgrams}
P_+ -A^*P_+A=C^*C \ands 
P_- -\alpha^*P_-\alpha =\gamma^*\gamma.
\end{equation} 
Note that in general the second and third term in the  realization  \eqref{R2} are not  required to be minimal realizations.
The next lemma provides the realization of  $R^*(z)R(z)$ starting from  \eqref{R2}.  

\begin{lemma}\label{lemRstarR}   The product function  $R^*(z)R(z)$ is given by 
\[
R^*(z)R(z)=\wt{D} + \wt{C}(z I-\wt{A})^{-1}\wt{B}+
z \wt{B}^*(I-z\wt{A}^*)^{-1} \wt{C}^*, 
\]
where 
\begin{align}
\wt{A}&=\begin{bmatrix}A^*& C^*\gamma\\0&\alpha \end{bmatrix}, \qquad\qquad\quad
\wt{B}=\begin{bmatrix} C^*R_0+A^*P_+B\\ \beta\end{bmatrix},\label{def:ABtilde}
\\[1mm]
\wt{C}&=\begin{bmatrix}B^*& R_0^*\gamma+\beta^*P_-\alpha\end{bmatrix},\quad \wt{D}= R_0^*R_0+B^*P_+B+\beta^*P_- \beta. \label{def:CDtilde}
\end{align}
Furthermore, the matrix $\wt{A}$ is stable.
\end{lemma}

\noindent
{\bf Proof.} We have 
$$
R^*(z)=R(1/ \bar{z})^*=R_0^*+B^*(zI-A^*)^{-1}C^*+z \beta^* (I-z\alpha^*)^{-1}\gamma^*,
$$
and so
\begin{align*} 
& R^*(z)R(z)=R_0^*R_0+zR_0^*C(I-zA)^{-1}B+R_0^*\gamma (zI-\alpha)^{-1}\beta +
\nonumber
\\ 
& +B^*(zI-A^*)^{-1}C^*R_0 +z \beta^* (I-z\alpha^*)^{-1}\gamma^*R_0 +
\nonumber
\\
&+zB^*(zI-A^*)^{-1}C^*C(I-zA)^{-1}B+z \beta^* (I-z\alpha^*)^{-1}\gamma^* \gamma (zI-\alpha)^{-1}\beta +
\nonumber
\\
&+B^*(zI-A^*)^{-1}C^*\gamma (zI-\alpha)^{-1}\beta + z^2\beta^* (I-z\alpha^*)^{-1}\gamma^*C(I-zA)^{-1}B.
\end{align*}

Let $P_+$ and $P_-$ be the observability gramians given by
\eqref{eq:obsgrams}.
Then (compare the argument following formula (2.17) in \cite{FKR})
\begin{align*}
&zB^*(zI-A^*)^{-1}C^*C(I-zA)^{-1}B=
\\
&zB^*P_+A(I-zA)^{-1}B + B^*P_+B +B^*(zI-A^*)^{-1}A^*P_+B,
\end{align*}
and
\begin{align*}
&z\beta^* (I-z\alpha^*)^{-1}\gamma^* \gamma (zI-\alpha)^{-1}\beta
= 
\\
&\beta^* P_-\alpha  (zI-\alpha)^{-1}\beta +\beta^* P_-\beta + z\beta^* (I-z\alpha^*)^{-1}
\alpha^* P_-\beta.
\end{align*}
Inserting this in the formula for $R^*(z)R(z)$ and regrouping terms a bit we obtain
\begin{align*}
& R^*(z)R(z)=R_0^*R_0+B^*P_+B+\beta^*P_-\beta +
\\
&+z(R_0^*C+B^*P_+A)(I-zA)^{-1}B+B^*(z I-A^*)^{-1}(C^*R_0+A^*P_+B) +
\\
&+(R_0^*\gamma +\beta^*P_-\alpha)(z I-\alpha)^{-1}\beta +
z\beta^*(I-z\alpha^*)^{-1}(\gamma^*R_0+\alpha^*P_-\beta) +
\\
&+B^*(zI-A^*)^{-1}C^*\gamma (zI-\alpha)^{-1}\beta + z^2\beta^* (I-z\alpha^*)^{-1}\gamma^*C(I-zA)^{-1}B.
\end{align*}
Now note that
\begin{align*}
&\begin{bmatrix}
B^* & R_0^*\gamma +\beta^*P_-\alpha 
\end{bmatrix}
\left( zI -\begin{bmatrix}
A^* & C^*\gamma \\
0 & \alpha
\end{bmatrix}\right)^{-1}
\begin{bmatrix}
C^*R_0 +A^*P_+B \\
\beta
\end{bmatrix}
\\
&=B^*(z I-A^*)^{-1}(C^*R_0 +A^*P_+B)+ (R_0^*\gamma +\beta^*P_-\alpha )(z I-\alpha)^{-1}\beta +\\
&+B^*(z I-A^*)^{-1}C^*\gamma(z I-\alpha)^{-1}\beta ,
\end{align*}
and
\begin{align*}
&z\begin{bmatrix}
R_0^*C+B^*P_+A & \beta^*
\end{bmatrix}
\left( I-z\begin{bmatrix}
A & 0 \\
\gamma^*C & \alpha^*
\end{bmatrix}\right)^{-1}
\begin{bmatrix}
B \\ 
\gamma^*R_0 + \alpha^*P_-\beta
\end{bmatrix}
\\
&=z(R_0^*C+B^*P_+A)(I-zA)^{-1}B + z\beta^*(I-z\alpha^*)^{-1}(\gamma^*R_0 + \alpha^*P_-\beta) +
\\
&+z^2 \beta^*(I-z\alpha^*)^{-1}\gamma^*C(I-zA)^{-1}B .
\end{align*}
So we arrive at the following formula for $R^*(z)R(z)$:
\begin{align}\label{eq:RstarR}
&R^*(z)R(z) = R_0^*R_0+B^*P_+B+\beta^*P_-\beta +\nonumber \\
& + \begin{bmatrix}
B^* & R_0^*\gamma +\beta^*P_-\alpha 
\end{bmatrix}
\left( z I -\begin{bmatrix}
A^* & C^*\gamma \\
0 & \alpha
\end{bmatrix}\right)^{-1}
\begin{bmatrix}
C^*R_0 +A^*P_+B \\
\beta
\end{bmatrix} +\nonumber
\\&+z\begin{bmatrix}
R_0^*C+B^*P_+A & \beta^*
\end{bmatrix}
\left( I-z\begin{bmatrix}
A & 0 \\
\gamma^*C & \alpha^*
\end{bmatrix}\right)^{-1}
\begin{bmatrix}
B \\ 
\gamma^*R_0 + \alpha^*P_-\beta
\end{bmatrix}.
\end{align}
This completes the proof of the formula for $R^*(z)R(z)$ upon using the definitions of $\wt{A}, \wt{B}, \wt{C}$, and $\wt{D}$.

The fact that $\wt{A}$ is stable follows from the stability of $A$ and $\alpha$.
 \hfill$\Box$
 
 \medskip
Since $R(z)$  is an $m\ts m$ rational matrix function with  no poles and no zeros on the unit circle,  and since  $A$ and $\alpha$ are stable, we can use Theorem 1.1 in \cite{FKR} to prove the  following lemma. In fact,  the following lemma follows from the symmetric version of   Theorem 1.1 in \cite{FKR}.

\begin{lemma}\label{lem:2.2}
 Let $\wt{A}$, $\wt{B}$,  $\wt{C}$, and $\wt{D}$  be the matrices defined by \eqref{def:ABtilde} and \eqref{def:CDtilde}.  Then the algebraic  Riccati equation 
\begin{equation}\label{eq:Ricc}
  Q-\wt{A}Q\wt{A}^* = \left(\wt{B}-\wt{A}Q\wt{C}^*\right)
  \left(\wt{D}-\wt{C} Q\wt{C}^*\right)^{-1}\left(\wt{B}^*-\wt{C}Q\wt{A}^*\right)
 \end{equation}
has a (unique) stabilizing solution $Q$.
\end{lemma}

\medskip
Given   the solution $Q$ of  \eqref{eq:Ricc} we define
\begin{align}
D_{\phantom{0}}&=\left(\wt{D}-\wt{C}Q\wt{C}^*\right)^{1/2}, \quad 
C_0=\wt{B}^*- \wt{C}Q\wt{A}^*\nonumber\\[1mm]
B_0&=\wt{C}^*, \qquad \qquad \qquad\quad
A_0=\wt{A}^* -  \wt{C}^* D^{-2}C_0.\label{A0B0C0}
\end{align}
Now we can apply Lemma \ref{lem:2.2} and the symmetric version of Theorem 1.1 in \cite{FKR} to obtain the following proposition.

\begin{proposition} \label{prop:realpsi}
The matrix $A_0$ is stable and $R^*(z)R(z)=\Psi^*(z)\Psi(z)$. Here 
$\Psi(z)$ is the outer function given by 
\begin{equation}\label{def:psi}
\Psi(z)= D+zD^{-1}C_0 \left(I-z\wt{A}^* \right)^{-1}B_0, 
\end{equation}
and the inverse of $\Psi(z)$ is given by
\begin{equation}\label{def:psi-inv}
\Psi(z)^{-1}=D^{-1}-zD^{-2}C_0(I-zA_0)^{-1}
B_0D^{-1}.
\end{equation}
\end{proposition}

\medskip
Observe  that the Riccati equation can be rewritten as a Stein equation, namely
\begin{align}
Q&=\begin{bmatrix}
A^* & C^*\gamma \\ 0 & \alpha
\end{bmatrix}QA_0
+\begin{bmatrix}  
C^*R_0+A^*P_+B \\ \beta
\end{bmatrix}D^{-2}C_0 \label{altRicc}\\[1mm]
&=\wt{A}QA_0+\wt{B}D^{-2}C_0. \nonumber
\end{align}

\setcounter{equation}{0}
\section{The function $\Xi(z)$  and the solution of the 
 Wiener-Hopf indices problem} \label{sec:three}

The function $\Xi(z)$ is defined by the formula
\begin{equation}\label{def:Xi}
\Xi(z):=R(z)\Psi(z)^{-1}
\end{equation}
where  $\Psi(z)$ is the invertible  outer function given by \eqref{def:psi}. Obviously, $\Xi(z)$ is a  rational $m\ts m$  matrix function.

\begin{lemma}\label{lem:Xi} The rational function $\Xi(z)$ is unitary for each $z\in \BT$, and the functions $R(z)$ and $\Xi(z)$  have the same right Wiener-Hopf  indices.\\
\end{lemma} 

\noindent
{\bf Proof.} The fact that $\Xi(z)$ takes unitary values on the unit circle is checked by direct computation. For $z\in \mathbb{T}$:
$$
\Xi(z)^*\Xi(z)=\Psi(z)^{-*}R(z)^*R(z)\Psi(z)^{-1}=I.
$$ 

Suppose that $R(z)=W_-(z)D(z)W_+(z)$ is a right Wiener-Hopf factorization of $R(z)$, then 
$\Xi(z)=W_-(z)D(z)\tilde{W}_+(z)$, where $\tilde{W}_+(z)=W_+(z)\Psi^{-1}(z)$. Since $\Psi(z)$ is an invertible outer function, $\tilde{W}_+$ has no poles and no zeros on the closed unit disc. Hence the factorization $\Xi(z)=W_-(z)D(z)\tilde{W}_+(z)$ is a right Wiener-Hopf factorization of $\Xi(z)$, and as a consequence $R(z)$ and $\Xi(z)$ have the same right Wiener-Hopf indices. \hfill$\Box$

\begin{definition} With some ambiguity (because $\Xi(z)$ is unique only up to multiplication on the right by a constant unitary matrix), the function $\Xi(z)$ is said to  be the  \emph{left unitary factor} of  the rational $m\ts m$ matrix function $R(z)$. 
\end{definition}

Our aim is to obtain the right Wiener-Hopf  indices of $R(z)$. Recall  that  the above lemma tells us that  $R(z)$ and $\Xi(z)$  have the same right Wiener-Hopf  indices. Thus  it suffices to obtain  the right Wiener-Hopf  indices of the left  unitary factor  $\Xi(z)$.  Since $\Xi(z)$ is unitary for each $z\in \BT$,  we can apply Theorem 1.1 of  \cite{GKR} to get these indices.

As a first step (see \cite[page 696]{GKR}) we use the fact that $\Xi(z)$ factors as   
$\Xi(z)=V(z)W^*(z)$,  where $V(z)$ and $W(z)$ are rational bi-inner $m\ts m$ matrix functions. Furthermore, we may assume (without loss of generality) that  
\begin{align}
V(z)&=D_V+zC_V(I_{n_-}-zA_V)^{-1}B_V, \label{realiV}\\ 
W(z)&=D_W+zC_W(I_{n_+}-zA_W)^{-1}B_W, \label{realiW}
\end{align}
with $A_V$ and $A_W$ being stable, and  with both realizations being  unitary. The latter means that   
\[ 
\begin{bmatrix} A_V & B_V \\ C_V & D_V\end{bmatrix} \ands
\begin{bmatrix} A_W & B_W \\ C_W & D_W\end{bmatrix}
\]
are   both unitary matrices.  These realizations are then minimal. In what follows  we also need the unique solution $X$ of   the Stein equation 
\begin{equation}\label{Stein}
X- A_VXA_W^*=B_VB_W^*, \quad \mbox{that is,} \quad  X=\sum_{j=0}^{\iy} A_V^jB_VB_W^* (A_W^j)^*
\end{equation}
(see \cite{FK} where this matrix $X$ was introduced).

Next, we need the following standard notation: for a pair $(C,A)$ we denote
${\rm Ker\, }_k(C , A)=\cap_{j=0}^{k-1} {\rm Ker\, }CA^j$, and for a pair $(A,B)$ we denote
${\rm Im\, }_k(A, B)\allowbreak = \vee_{j=0}^{k-1} {\rm Im\,}A^jB$. For $k=0$ these spaces are interpreted as the full Euclidean space, respectively the zero subspace.  

By a direct application of Theorem 1.1 in \cite{GKR} we have the following result on the Wiener-Hopf indices of $R(z)$, as they are the same as the Wiener-Hopf indices of $\Xi(z)$.

\begin{theorem}\label{thm:main} 
With notation as above the number $s$ of negative right Wiener-Hopf indices of  the  left  unitary factor  $\Xi(z)$ of $R(z)$ is given by
\[
s={\rm dim\,}{\rm Ker\, }X-{\rm dim\,} {\rm Ker\,}\begin{bmatrix} B_W^*\\ XA_W^* \end{bmatrix},
\]
while the negative right Wiener-Hopf indices $-\alpha_j$, $j=1,\ldots ,s$, are given by
\begin{align*}
&\alpha_j=\# \Big{\{} k\in\mathbb{N} \mid
{\rm dim\ } {\rm Ker}_{k-1}\left(\begin{bmatrix} B_W^*\\ XA_W^*\end{bmatrix}, A_W^*\right) - \\
&\hspace{4cm}-{\rm dim\,} {\rm Ker}_{k}\left(\begin{bmatrix} B_W^*\\ XA_W^* \end{bmatrix} , A_W^*\right) \geq j \Big{\}}.
\end{align*} 
Furthermore, the number $t$ of positive right Wiener-Hopf indices of  $\Xi(z)$ is given by
\[
t={\rm dim\,}{\rm Im\, }\begin{bmatrix}  B_V & A_V X \end{bmatrix}- {\rm dim\, }{\rm Im\, }X
\]
while the positive right Wiener-Hopf indices $\omega_j$, $j=1, \ldots , t$, are given by
\begin{align*}
&\omega_j =\# \Big{\{} k\in\mathbb{N} \mid \dim{\rm Im\,}_k\left(A_V ,
\begin{bmatrix} B_V &A_V X\end{bmatrix}\right)-    \\
&\hspace{3cm}-\dim
{\rm Im\,}_{k-1}\left(A_V , \begin{bmatrix} B_V &A_V X\end{bmatrix}\right) \geq j\Big{\}} .
\end{align*}
\end{theorem}

To make the Wiener-Hopf  indices in the above theorem more  explicit  we need formulas for the matrices $A_V$, $A_W$,   $B_V$, $B_W$, and $X$ in terms of the matrices occuring in the realization \eqref{R} of $R(z)$.
This will be done in the next  two sections.

\setcounter{equation}{0}
\section{Properties of the left unitary  factor  $\Xi(z)$}
In this section we first produce a realization  for the   unitary factor function $\Xi(z)$. Recall that $\Xi(z)=R(z)\Psi(z)^{-1}$, which has unitary values on the unit circle, i.e., $\Xi(z)^{-1}=\Xi(z)^*$ for $|z|=1$.
Let $Y$ be the solution to the Stein equation
\begin{equation}\label{eq:SteinY}
Y-\alpha YA_0=\alpha\beta D^{-2}C_0.
\end{equation}
Introduce the following notation:
\begin{align}
\Xi_0&=R_0D^{-1}-\gamma YA_0 B_0D^{-1}
-\gamma\beta D^{-2}C_0B_0D^{-1},\nonumber\\
C_1&=\left(\begin{bmatrix} 
C & 0 \end{bmatrix} -R_0D^{-2}C_0 - \gamma Y A_0^2-\gamma\beta D^{-2}C_0A_0
\right),\nonumber\\
\beta_1&=\left(\beta - Y B_0\right)D^{-1}.\label{Xi0C1beta1}
\end{align}
Observe that $\Xi_0, C_1$ and $\beta_1$ are defined in terms of operators which we know from the preceding sections and the solution $Y$ of the Stein equation \eqref{eq:SteinY}. Note also that the operators from the previous section in turn are explicitly expressed in terms of the operators in the realization of $R(z)$ and solutions to Stein equations and a Riccati equation, which can be computed.
With this notation we can give an explicit formula for the unitary function $\Xi(z)$ such that $R(z)=\Xi(z)\Psi(z)$.

\begin{proposition}
The left unitary factor  $\Xi(z)$ of  the rational $m\ts m$ matrix function $R(z)$ is given in realized form by the following formula:
\begin{equation}\label{realizationXi}
\Xi(z)=\Xi_0+\gamma(z I-\alpha)^{-1}\beta_1+zC_1\left( I-z  A_0 \right)^{-1}
 B_0D^{-1} .
\end{equation}
\end{proposition}

\noindent
{\bf Proof.} Compute
\begin{align*}
\Xi(z)=&R(z)\Psi(z)^{-1} = R_0D^{-1}+zC(I-zA)^{-1}BD^{-1} +\gamma(z I-\alpha)^{-1}\beta D^{-1}
\\
&-zR_0D^{-2}C_0(I-zA_0)^{-1}\begin{bmatrix}
B \\ 
\gamma^*R_0 + \alpha^*P_-\beta 
\end{bmatrix}D^{-1} \\
&-z\gamma(z I-\alpha)^{-1}\beta D^{-2}C_0(I-zA_0)^{-1}
\begin{bmatrix}
B \\ 
\gamma^*R_0 + \alpha^*P_-\beta 
\end{bmatrix}D^{-1} \\
&-z^2C(I-zA)^{-1}BD^{-2}C_0(I-zA_0)^{-1}
\begin{bmatrix}
B \\ 
\gamma^*R_0 + \alpha^*P_-\beta 
\end{bmatrix}D^{-1}.
\end{align*}
We can combine terms as follows:
\begin{align*}
&\Xi(z) = R_0D^{-1}+
z\begin{bmatrix} 
C & -R_0D^{-2}C_0
\end{bmatrix}
\left( I-z \begin{bmatrix} A & -BD^{-2}C_0 \\ 0 & A_0 \end{bmatrix}\right)^{-1}
\begin{bmatrix} B \\ B_0\end{bmatrix}D^{-1} \\
& + \gamma(z I-\alpha)^{-1}\beta D^{-1} -
z\gamma(z I-\alpha)^{-1}\beta D^{-2}C_0(I-zA_0)^{-1}B_0D^{-1}.
\end{align*}

Now from the Riccati equation \eqref{altRicc} (second row) we obtain 
\begin{align*}
\beta D^{-2}C_0 & =
\begin{bmatrix} Q_{21} & Q_{22}\end{bmatrix}
-\begin{bmatrix} 0 & \alpha \end{bmatrix} Q A_0
\end{align*}
which will be used in the last two terms of $\Xi(z)=R(z)\Psi(z)^{-1}$.
Indeed, one has
\begin{align*}
& z(z I-\alpha)^{-1}\beta D^{-2}C_0(I-zA_0)^{-1}= \\
= &z(z I-\alpha)^{-1}\Big(\begin{bmatrix} Q_{21} & Q_{22}\end{bmatrix}
-\begin{bmatrix} 0 & \alpha \end{bmatrix} Q A_0\Big)(I-zA_0)^{-1} \\
=&z(z I-\alpha)^{-1}\Big(\begin{bmatrix} Q_{21} & Q_{22}\end{bmatrix}
-\alpha \begin{bmatrix} Q_{21} & Q_{22} \end{bmatrix}  A_0\Big)(I-zA_0)^{-1} \\
= & (z I-\alpha)^{-1}\Big(z\begin{bmatrix} Q_{21} & Q_{22}\end{bmatrix}
-\alpha \begin{bmatrix} Q_{21} & Q_{22} \end{bmatrix} z A_0\Big)(I-zA_0)^{-1} \\
=& (z I-\alpha)^{-1}\Big(z\begin{bmatrix} Q_{21} & Q_{22}\end{bmatrix}-\alpha \begin{bmatrix} Q_{21} & Q_{22} \end{bmatrix} +\alpha \begin{bmatrix} Q_{21} & Q_{22} \end{bmatrix} -\\
&\hspace{5cm}-\alpha \begin{bmatrix} Q_{21} & Q_{22} \end{bmatrix} z A_0\Big)(I-zA_0)^{-1}  \\
=&  (z I-\alpha)^{-1}\Big((z I-\alpha)\begin{bmatrix} Q_{21} & Q_{22}\end{bmatrix}
+\alpha \begin{bmatrix} Q_{21} & Q_{22} \end{bmatrix}(I- z A_0)\Big)(I-zA_0)^{-1}\\
= &  \begin{bmatrix} Q_{21} & Q_{22} \end{bmatrix}(I- z A_0)^{-1} + (z I-\alpha)^{-1}\alpha \begin{bmatrix} Q_{21} & Q_{22} \end{bmatrix} \\
= & \begin{bmatrix} Q_{21} & Q_{22} \end{bmatrix}(I- z A_0)^{-1} + (z I-\alpha)^{-1} \begin{bmatrix}0 & \alpha \end{bmatrix}Q.
\end{align*}
Hence $\Xi(z)$ is given by
\begin{align*}
\Xi(z)= R_0 D^{-1}&+
z\begin{bmatrix} 
C &  -R_0D^{-2}C_0
\end{bmatrix}
\left( I-z \begin{bmatrix} A & -BD^{-2}C_0 \\ 0 & A_0 \end{bmatrix}\right)^{-1}
\begin{bmatrix} B \\ B_0\end{bmatrix}D^{-1} \\
& + \gamma(z I-\alpha)^{-1}\beta D^{-1}\\
&- \gamma(z I-\alpha)^{-1}\begin{bmatrix}0 & \alpha \end{bmatrix}QB_0D^{-1}\\
&-\gamma \begin{bmatrix} Q_{21} & Q_{22} \end{bmatrix}(I- z A_0)^{-1} B_0D^{-1}.
\end{align*}
Rewrite the last term as follows:
\begin{align*}
\gamma \begin{bmatrix} Q_{21} & Q_{22} \end{bmatrix}&(I- z A_0)^{-1} B_0D^{-1}=\\
=&\begin{bmatrix}0 & \gamma \end{bmatrix}Q(I- z A_0)^{-1} B_0D^{-1}\\
=&\begin{bmatrix}0 & \gamma \end{bmatrix}QB_0D^{-1}+z\begin{bmatrix}0 & \gamma \end{bmatrix}QA_0(I- z A_0)^{-1} B_0D^{-1}.
\end{align*}
Then we arrive at the following formula for $\Xi(z)$:
\begin{align*}
\Xi(z)&=R_0D^{-1}-\begin{bmatrix}0 & \gamma \end{bmatrix}QB_0D^{-1}+
\\
&\ \ +\gamma(z I-\alpha)^{-1}\left(\beta -\begin{bmatrix}0 & \alpha \end{bmatrix}QB_0\right)D^{-1}\\
&\ \ +z\begin{bmatrix} 
C & -R_0D^{-2}C_0
\end{bmatrix}
\left( I-z \begin{bmatrix} A &- BD^{-2}C_0 \\ 0 & A_0 \end{bmatrix}\right)^{-1}
\begin{bmatrix} B \\ B_0\end{bmatrix}D^{-1} \\
&\ \  -z\begin{bmatrix}0 & \gamma \end{bmatrix}QA_0(I- z A_0)^{-1} B_0D^{-1}\\
&= R_0D^{-1}-\begin{bmatrix}0 & \gamma \end{bmatrix}QB_0D^{-1}
\\
&\ \ +\gamma(z I-\alpha)^{-1}\left(\beta -\begin{bmatrix}0 & \alpha \end{bmatrix}QB_0\right)D^{-1}\\
&\ \ +z\begin{bmatrix} 
C & \left(-R_0D^{-2}C_0 -\begin{bmatrix}0 & \gamma \end{bmatrix}QA_0\right)
\end{bmatrix}\ts\\
&\hspace{3cm} \ts \left( I-z \begin{bmatrix} A & -BD^{-2}C_0 \\ 0 & A_0 \end{bmatrix}\right)^{-1}
\begin{bmatrix} B \\ B_0\end{bmatrix}D^{-1} .
\end{align*}

Next, consider the last term in this expression, i.e.
$$
z\begin{bmatrix} 
C & \left(-R_0D^{-2}C_0 -\begin{bmatrix}0 & \gamma \end{bmatrix}QA_0\right)
\end{bmatrix}
\left( I-z \begin{bmatrix} A & -BD^{-2}C_0 \\ 0 & A_0 \end{bmatrix}\right)^{-1}
\begin{bmatrix} B \\ B_0\end{bmatrix}D^{-1} .
$$
Let $\Pi=\begin{bmatrix} I & 0\end{bmatrix}$, and let $S=\begin{bmatrix} I & \Pi \\ 0 & I \end{bmatrix}$. Observe that $-BD^{-2}C_0 -\Pi A_0 +A\Pi =0$. Hence
$$
S^{-1}\begin{bmatrix} A & -BD^{-2}C_0 \\ 0 & A_0 \end{bmatrix}S =
\begin{bmatrix} A & 0 \\ 0 & A_0 \end{bmatrix}, \qquad S^{-1}\begin{bmatrix} B\\ B_0 \end{bmatrix}=\begin{bmatrix} 0 \\ B_0 \end{bmatrix}.
$$
Hence, applying a similarity transformation with $S$ on the realization of the last term we arrive at
\begin{align*}
\Xi(z)&= R_0D^{-1}-\begin{bmatrix}0 & \gamma \end{bmatrix}QB_0D^{-1}
\\
&\ \ +\gamma(z I-\alpha)^{-1}\left(\beta -\begin{bmatrix}0 & \alpha \end{bmatrix}QB_0\right)D^{-1}\\
&\ \ +z\left(\begin{bmatrix} 
C & 0 \end{bmatrix} -R_0D^{-2}C_0 -\begin{bmatrix}0 & \gamma \end{bmatrix}QA_0
\right)
\left( I-z  A_0 \right)^{-1}
 B_0D^{-1} .
\end{align*}

To rewrite the formula once more, consider the form \eqref{altRicc} of the algebraic Riccati equation
and premultiply by $\begin{bmatrix}0 & \alpha \end{bmatrix}$:
$$
\begin{bmatrix}0 & \alpha \end{bmatrix}Q=\begin{bmatrix}0 & \alpha^2 \end{bmatrix}QA_0+\alpha\beta D^{-2}C_0.
$$
Repeating the argument, premultiply by $\begin{bmatrix}0 & \alpha^2 \end{bmatrix}$ and insert the resulting expression to obtain
$$
\begin{bmatrix}0 & \alpha \end{bmatrix}Q=\begin{bmatrix}0 & \alpha^3 \end{bmatrix}QA_0^2+\alpha^2\beta D^{-2}C_0A_0+\alpha\beta D^{-2}C_0.
$$
Continuing by induction we have 
$$
\begin{bmatrix}0 & \alpha \end{bmatrix}Q=\begin{bmatrix}0 & \alpha^{n+1} \end{bmatrix}QA_0^n+
\sum_{j=1}^n \alpha^j\beta D^{-2}C_0A_0^{j-1}.
$$
Since $\alpha$ and $A_0$ are stable we have $\lim_{n\to\infty}\alpha^n=0$ and $\lim_{n\to\infty}A_0^n=0$, and so
$$
\begin{bmatrix}0 & \alpha \end{bmatrix}Q=\sum_{j=1}^\infty \alpha^j\beta D^{-2}C_0a_0^{j-1}=
Y.
$$
In addition, again using \eqref{altRicc} yields
$$
\begin{bmatrix}0 & \gamma \end{bmatrix}Q=\gamma\left(\begin{bmatrix}0 & \alpha \end{bmatrix}QA_0+\beta D^{-2}C_0\right)=\gamma YA_0+\gamma\beta D^{-2}C_0.
$$

Inserting this in the formula for $\Xi(z)$ we have
\begin{align*}
\Xi(z)&= R_0D^{-1}-\gamma YA_0B_0D^{-1}
-\gamma\beta D^{-2}C_0B_0D^{-1} \\
&\quad +\gamma(z I-\alpha)^{-1}\left(\beta -YB_0\right)D^{-1}\\
&\quad+z\Big(\begin{bmatrix} C & 0 \end{bmatrix} -R_0D^{-2}C_0 - \gamma YA_0^2
 - \gamma\beta D^{-2}C_0A_0 \Big)
\left( I-z  A_0 \right)^{-1}B_0D^{-1} .
\end{align*}
Using the notations introduced in the beginning of this section, the formula for $\Xi(z)$ can be rewritten as in the statement of the proposition. \hfill $\Box$

\medskip
{\bf Remark.}
Instead of \eqref{realizationXi} the realization of  $\Xi(z)$  can also be obtained  by using the matrices $A_V$ and $A_W^*$  in \eqref{realiV} and  \eqref{realiW} and the matrix $X$ in \eqref{Stein} in place of the matrices  $A_0$ and $\alpha$.  More precisely, see, e.g., \cite{FK},
\begin{equation}\label{realizationXi2}
\Xi(z)={\Xi}_0 + C_R(zI -A_W^*)^{-1}C_W^*+ zC_V(I-zA_V)^{-1}B_R ,
\end{equation}
where
\begin{align*}
B_R& = B_VD_W^*+A_VXC_W^* ,\\
C_R&=D_VB_W^* +C_VXA_W^*,\\
{\Xi}_0&=D_VD_W^*+C_VXC_W^*.  \qquad \Box
\end{align*}

\medskip
Using the fact that $\Xi(z)$ is a unitary function on $\BT$, we have the following extra properties for the matrices in the realization, some of which require the extra condition that $A_0^*$ and $\alpha$ have no common eigenvalues, and that, in addition, the pair $(\alpha, \beta_1)$ is controllable.

\begin{proposition}\label{propertiesunitary}
Let $\Xi(z)$ be the function \eqref{realizationXi} taking unitary values on the unit circle. Assume that the pair $(\gamma,\alpha)$ is observable, and the pair $(A_0,B_0)$ is controllable. Then the unsymmetric Lyapunov equation
\begin{equation}\label{Lyap}
A_0^*P_1-P_1\alpha=-C_1^*\gamma,
\end{equation}
has a solution. One particular solution is the zero-pole coupling matrix of $\Xi(z)$ corresponding to the pole pair $(\gamma, \alpha)$ and the zero pair $(A_0^*,C_1^*)$ corresponding to the unit disc; this solution will be denoted by $P_1$ henceforth.

Let $P_0$ be the solution to the Stein equation
\begin{equation}\label{SteinP0}
P_0-A_0^*P_0A_0=C_1^*C_1.
\end{equation}

Assume further that $A_0^*$ and $\alpha$ have no common eigenvalues. Then \eqref{Lyap} has a unique solution, and the following identities hold
\begin{eqnarray}
& \Xi_0^*\Xi_0+D^{-*}B_0^*P_0B_0D^{-1}+\beta_1^*P_-\beta_1 & = I, \label{unitarycond1}\\
&C_1^*\Xi_0+A_0^*P_0B_0D^{-1}-P_1\beta_1 & = 0. \label{unitarycond2}
\end{eqnarray}
If in addition the pair $(\alpha,\beta_1)$ is controllable, then also
\begin{equation}
D^*\Xi_0^*\gamma+D^*\beta_1^*P_-\alpha +B_0^*P_1  =0.\label{unitarycond3}
\end{equation}

\end{proposition}

\noindent
{\bf Proof.} Since $(\gamma,\alpha)$ is observable and $\alpha$ and $A_0$ are stable matrices, a pole pair of $\Xi(z)$ corresponding to the unit disc is given by $(\gamma, \alpha)$. Likewise, considering $\Xi^*(z)=\Xi(\tfrac{1}{\bar{z}})^*$, we have
$$
\Xi^*(z)=\Xi_0^*+D^{-*}B_0^*(z I-A_0^*)^{-1}C_1^*+z\beta_1^*(I-z\alpha^*)^{-1}\gamma^*,
$$
and since $(A_0,B_0)$ is controllable and $\Xi^*(z)=\Xi(z)^{-1}$, a zero pair of $\Xi(z)$ corresponding to the unit disc is given by $(A_0^*,C_1^*)$. Let $P_1$ be the zero-pole coupling matrix, then $P_1$ satisfies \eqref{Lyap} (see, \cite{BGR}, Sections 4.5 and 4.4, see also \cite{Ball-Ran1, Ball-Ran2}).

Apply Lemma \ref{lemRstarR} to $\Xi(z)$ to obtain
\begin{align*}
\Xi^*(z)\Xi(z)  = & \Xi_0^*\Xi_0+D^{-*}B_0^*P_0B_0D^{-1}+\beta_1^*P_-\beta_1 \\
& +\widehat{C}(z I-\widehat{A})^{-1}\widehat{B}+z\widehat{B}^*(I-z\widehat{A}^*)^{-1}\widehat{C}^*,
\end{align*}
where
\begin{align*}
\widehat{C} & = \begin{bmatrix} D^{-*}B_0^* & \Xi_0^*\gamma + \beta_1^*P_-\alpha\end{bmatrix} ,\\
\widehat{A}& = \begin{bmatrix} A_0^* & C_1^*\gamma \\ 0 & \alpha \end{bmatrix}, \qquad
\widehat{B} = \begin{bmatrix} C_1^*\Xi_0+A_0^*P_0B_0D^{-1} \\ \beta_1\end{bmatrix}.
\end{align*}
Since $\Xi(z)$ is unitary on $\BT$ it is now immediate that \eqref{unitarycond1} holds, and moreover, that
\begin{align*}
0  = &\widehat{C}(z I-\widehat{A})^{-1}\widehat{B} \\
= &D^{-*}B_0^*(z I-A_0^*)^{-1}(C_1^*\Xi_0+A_0^*P_0B_0D^{-1} )\\
&+ (\Xi_0^*\gamma + \beta_1^*P_-\alpha)(z I-\alpha)^{-1}\beta_1 \\
&+D^{-*}B_0^*(z I-A_0^*)^{-1}C_1^*\gamma(z I-\alpha)^{-1}\beta_1.
\end{align*}

Let $P_1$ be the zero-pole coupling matrix, which is a solution to \eqref{Lyap}. Then
$$
(z I-A_0^*)^{-1}C_1^*\gamma(z I-\alpha)^{-1} =- (z I-A_0^*)^{-1}P_1+P_1(z I-\alpha)^{-1}
$$
and inserting this into the formula for $\widehat{C}(z I-\widehat{A})^{-1}\widehat{B}$ we obtain
\begin{align*}
0= &D^{-*}B_0^*(z I-A_0^*)^{-1}(C_1^*\Xi_0+A_0^*P_0B_0D^{-1}-P_1\beta_1 )\\
&+ (\Xi_0^*\gamma + \beta_1^*P_-\alpha+D^{-*}B_0^*P_1)(z I-\alpha)^{-1}\beta_1.
\end{align*}

Now using the assumption that
$A_0^*$ and $\alpha$ have no common eigenvalues each of the two terms must be zero. 
By controllability of $(A_0,B_0)$ \eqref{unitarycond2} follows, and under the assumption that $(\alpha, \beta_1)$ is controllable also \eqref{unitarycond3} follows.
\hfill $\Box$

\medskip
We conclude this section with the special case when 
\begin{equation}\label{specialR}
R(z)=R_0+zC(I-zA)^{-1}B\quad\mbox{with $A$ being stable}. 
\end{equation}
This yields the following corollary.

\begin{corollary} Assume $R(z)$ is given by \eqref{specialR}. Then the left unitary factor $\Xi(z)$ of $R(z)$ is a bi-inner function, and is given by
\[
\Xi(z)=R_0D^{-1}+z(C-R_0D^{-2}C_0)(I-zA_0)^{-1}BD^{-1}.
\]
\end{corollary}

\noindent
{\bf Proof.}
In the special case where $R(z)$ is stable, the algebraic Riccati equation \eqref{eq:Ricc} becomes
$$
Q=A^*QA+(C^*R_0+A^*(P_+-Q)B)(R_0^*R_0+B^*(P_+-Q)B)^{-1}(R_0^*C+B^*(P_+-Q)A)
$$
and for the stabilizing solution $Q$ we have
\begin{align*}
D&= (R_0^*R_0+B^*(P_+-Q)B)^{1/2},\\
C_0&=R_0^*C+B^*(P_+-Q)A,\\
A_0&=A-BD^{-2}C_0.
\end{align*}
Moreover, $B_0=B$.

Thus the outer factor $\Psi(z)$ and its inverse are given by
\begin{align*}
\Psi(z)&=D+zD^{-1}C_0(I-zA)^{-1}B, \\
\Psi(z)^{-1}&=D^{-1} -zD^{-2}C_0(I-zA_0)^{-1}BD^{-1},
\end{align*}
while the unitary factor is given by
$$
\Xi(z)=R_0D^{-1}+z(C-R_0D^{-2}C_0)(I-zA_0)^{-1}BD^{-1}.
$$
Note that $\Xi(z)$ is stable, so this is an inner factor, as expected for the case where $R(z)$ is stable.
In fact, because $R(z)$ is assumed to be square, $\Xi(z)$ is bi-inner.
Compare \cite{FB} and \cite{FR}. \hfill$\Box$

\setcounter{equation}{0}
\section{Construction of the factorization $\Xi(z)=V(z)W^*(z)$}

In this section $\Xi(z)$ is the left unitary factor of $R(z)$,  and $V(z)$ and $W(z)$  are the  bi-inner $m\ts m$ matrix functions appearing in \eqref{realiV} and \eqref{realiW}.  The identity  $\Xi(z)=V(z)W^*(z)$  is a Douglas-Shapiro-Shields  factorization of $\Xi(z)$ (see  Section 4.7 in \cite{FB}; see also \cite{Potapov-60}).  In this section we shall construct  explicit formulas for the matrix functions  $V(z)$ and $W(z)$. 

\begin{proposition}\label{DSSexplicit}
Let $\Xi(z)$ be given by \eqref{realizationXi}. Then the bi-inner rational matrix functions $V(z)$ and $W(z)$ in the   factorization $\Xi(z)=V(z)W^*(z)$  can be constructed as follows.

Let  $P_0$ be the  solution of the Stein equation
\[
P_0-A_0^*P_0A_0=C_1^*C_1,\  \mbox{that is,\  $P_0=\sum_{j=0}^\iy (A_0^*)^j C_1^*C_1A_0^j.$}
\]
Then there exist matrices    $B_V$ and $D_V$ such that
\begin{equation}\label{uniV}
\begin{bmatrix}
A_0^* & C_1^* \\ B_V^* & D_V^*
\end{bmatrix}
\begin{bmatrix}
P_0 & 0 \\ 0 & I 
\end{bmatrix}
\begin{bmatrix}
A_0 & B_V \\ C_1 & D_V 
\end{bmatrix}
=
\begin{bmatrix}
P_0 & 0 \\ 0 & I
\end{bmatrix}.
\end{equation}
Moreover,  the function 
\begin{equation}\label{formulaV}
V(z)= D_V + zC_1(I-zA_0)^{-1}B_V
\end{equation}
is a bi-inner  factor  for $\Xi(z)$ and the function $W(z)$ defined  by $W^*(z)=V^*(z)\Xi(z)$ is  bi-inner too.

Assume, in addition,  that $A_0^*$  and $\alpha$ have no common eigenvalues and that the realization \eqref{realizationXi} is minimal. Let $P_1$ be the solution of \eqref{Lyap}.
Set
\[
D_W^* =  D_V^*\Xi_0 + B_V^*P_0B_0D^{-1} \ands 
B_W^*=D_V^*\gamma+B_V^*P_1.
\]
Then
\begin{equation}\label{formulaW}
W^*(z)= D_W^* 
+B_W^*(z I-\alpha)^{-1}\beta_1.
\end{equation}
\end{proposition}

\medskip
Note that finding $B_V$ and $D_V$ is a straightforward completion problem. Indeed, the columns of
$\begin{bmatrix} P_0^{1/2}A_0\\ C_1\end{bmatrix}$ are orthonormal, and the columns of $\begin{bmatrix} B_V\\ D_V\end{bmatrix}$ complete the set of columns of $\begin{bmatrix} P_0^{1/2}A_0\\ C_1\end{bmatrix}$ to an orthonormal basis.

\noindent
{\bf Proof.}
To construct the factors in the DSS factorization explicitly from a realization of $\Xi(z)$ we follow the procedure described in Theorem 4.2.1, Remark 4.3.4 and Section 4.8.1  in \cite{FB}. That procedure leads directly to the formula \eqref{formulaV} for $V(z)$. The formula may also be derived from the results in Chapter 7 of \cite{BGR}, assuming the pair $(C_1, A_0)$ is observable. Indeed, $V(z)$ has the pole pair $(C_1, A_0)$  corresponding to the unit disc, and has to be inner. Then $P_0$ has to be positive definite, and there are $B_V$ and $D_V$ satisfying \eqref{uniV}, while $V(z)$ is given by \eqref{formulaV}.

It remains to compute $W(z)$. This will be done by direct computation of the product $W^*(z)=V^*(z)\Xi(z)$.
Let us compute $V^*(z)\Xi(z)$ for the given $V(z)$:
\begin{align*}
W^*(z)&=V^*(z)\Xi(z)=(D_V^*+B_V^*(z I-A_0^*)^{-1}C_1^*)\times 
\\
&(\Xi_0+\gamma(z I-\alpha)^{-1}\beta_1+zC_1(I-zA_0)^{-1}B_0D^{-1}) \\
=& D_V^*\Xi_0 + D_V^*\gamma(z I-\alpha)^{-1}\beta_1+zD_V^*C_1(I-zA_0)^{-1}B_0D^{-1} \\
&+B_V^*(z I-A_0^*)^{-1}C_1^*\Xi_0+B_V^*(z I-A_0^*)^{-1}C_1^*\gamma(z I-\alpha)^{-1}\beta_1 \\
&+zB_V^*(z I-A_0^*)^{-1}C_1^*C_1(I-zA_0)^{-1}B_0D^{-1}.
\end{align*}
Now the terms $D_V^*\gamma(z I-\alpha)^{-1}\beta_1$ and 
$ B_V^*(z I-A_0^*)^{-1}C_1^*\Xi_0 $ as well as $B_V^*(z I-A_0^*)^{-1}C_1^*\gamma(z I-\alpha)^{-1}\beta_1$
have only Fourier coefficients corresponding to negative powers of $z$, while $zD_V^*C_1(I-zA_0)^{-1}B_0D^{-1}$ has only Fourier coefficients corresponding to positive powers of $z$. 
So the constant term of $V^*(z)\Xi(z)$ is equal to $ D_V^*\Xi_0 $ plus the constant term of 
$$zB_V^*(z I-A_0^*)^{-1}C_1^*C_1(I-zA_0)^{-1}B_0D^{-1}.$$

We re-express this term with the help of $P_0$, using
$$
(z I-A_0^*)^{-1}C_1^*C_1(I-zA_0)^{-1}=P_0A_0(I-zA_0)^{-1} +(z I-A_0^*)^{-1}P_0.
$$
this implies
\begin{align*}
&zB_V^*(z I-A_0^*)^{-1}C_1^*C_1(I-zA_0)^{-1}B_0D^{-1}
\\
=& zB_V^*P_0A_0(I-zA_0)^{-1}B_0D^{-1} +zB_V^*(z I-A_0^*)^{-1}P_0B_0D^{-1}.
\end{align*}
The first term in the latter expression again has only Fourier coefficients corresponding to positive powers of $z$. Furthermore,
\begin{align*}
zB_V^*(z I-A_0^*)^{-1}P_0B_0D^{-1}&=B_V^*(z I-A_0^*+A_0^*)(z I-A_0^*)^{-1}P_0B_0D^{-1}
\\ 
&=B_V^*P_0B_0D^{-1}+B_V^*(z I-A_0^*)^{-1}A_0^*P_0B_0D^{-1}.
\end{align*}
So the constant term of $V^*(z)\Xi(z)$ is equal to $ D_V^*\Xi_0 + B_V^*P_0B_0D^{-1}$.
It follows that 
$$
D_W^*=D_V^*\Xi_0 + B_V^*P_0B_0D^{-1}.
$$

We carry the computation a bit further to compute $W^*(z)$ directly, but it requires a bit of extra work to see that many terms cancel.

In fact, we obtain 
\begin{align*}
W^*(z)&= D_V^*\Xi_0 + B_V^*P_0B_0D^{-1} + D_V^*\gamma(z I-\alpha)^{-1}\beta_1+
\\
&+zD_V^*C_1(I-zA_0)^{-1}B_0D^{-1} \\
&+B_V^*(z I-A_0^*)^{-1}C_1^*\Xi_0\\
&+B_V^*(z I-A_0^*)^{-1}C_1^*\gamma(z I-\alpha)^{-1}\beta_1 \\
&+B_V^*(z I-A_0^*)^{-1}A_0^*P_0B_0D^{-1} \\
&+zB_V^*P_0A_0(I-zA_0)^{-1}B_0D^{-1}.
\end{align*}
As $B_V^*P_0A_0=-D_V^*C_1$ by \eqref{uniV}, the two terms involving $(I-zA_0)^{-1}$ cancel, and we are left with
\begin{align*}
W^*(z)&= D_V^*\Xi_0 + B_V^*P_0B_0D^{-1} + D_V^*\gamma(z I-\alpha)^{-1}\beta_1+
\\
&+B_V^*(z I-A_0^*)^{-1}C_1^*\Xi_0\\
&+B_V^*(z I-A_0^*)^{-1}C_1^*\gamma(z I-\alpha)^{-1}\beta_1 \\
&+B_V^*(z I-A_0^*)^{-1}A_0^*P_0B_0D^{-1}.
\end{align*}
Rewrite this by combining the fourth and sixth terms:
\begin{align*}
W^*(z)&= D_V^*\Xi_0 + B_V^*P_0B_0D^{-1} + D_V^*\gamma(z I-\alpha)^{-1}\beta_1+
\\
&+B_V^*(z I-A_0^*)^{-1}(C_1^*\Xi_0+A_0^*P_0B_0D^{-1})\\
&+B_V^*(z I-A_0^*)^{-1}C_1^*\gamma(z I-\alpha)^{-1}\beta_1.
\end{align*}

Next, we assume in addition that $(\gamma,\alpha)$ and $(A_0,B_0)$ are, respectively, observable and controllable, i.e., the realization \eqref{realizationXi} is minimal.
Let $P_1$ be the corresponding zero-pole coupling matrix.
Then $P_1$ is a solution of
$$
A_0^*P_1-P_1\alpha = - C_1^*\gamma.
$$
Hence, as in the proof of Proposition \ref{propertiesunitary}
\begin{align*}
W^*(z)&= D_V^*\Xi_0 + B_V^*P_0B_0D^{-1} \\
&+(D_V^*\gamma+B_V^*P_1)(z I-\alpha)^{-1}\beta_1 \\
&+B_V^*(z I-A_0^*)^{-1}(C_1^*\Xi_0+A_0^*P_0B_0D^{-1}-P_1\beta_1).
\end{align*}
Using the assumption that $A_0^*$ and $\alpha$ have no common eigenvalues we can use \eqref{unitarycond2}, and so we have that this equals
$$
W^*(z)= D_V^*\Xi_0 + B_V^*P_0B_0D^{-1} 
+(D_V^*\gamma+B_V^*P_1)(z I-\alpha)^{-1}\beta_1.
$$
Hence we have
\[
D_W^*  =  D_V^*\Xi_0 + B_V^*P_0B_0D^{-1}\ands
B_W^*=D_V^*\gamma+B_V^*P_1.
\]
This completes the proof.\hfill $\Box$

\medskip

In order to apply Theorem \ref{thm:main} we need the matrix $X$ in equation \eqref{Stein}, applied to the realizations \eqref{formulaV} and \eqref{formulaW}. The following proposition gives an explicit expression for this matrix.

\begin{proposition}\label{formulaX}
Let $V(z)$ and $W(z)$ be given by \eqref{formulaV} and \eqref{formulaW}, respectively. Then the corresponding matrix $X$ solving the equation 
\begin{equation}\label{equationX2}
X-A_0X\alpha=B_VB_W^*
\end{equation}
is given by 
\begin{equation}\label{eq:X2}
X=P_0^{-1}P_1.
\end{equation}
\end{proposition} 

\noindent
{\bf Proof.} Since $X$ is the unique solution of \eqref{equationX2}, all we need to do is to check that 
$P_0^{-1}P_1$ satisfies this equation. This is done by direct checking. First we note that by taking inverses in \eqref{uniV} and re-arranging terms
we have that
$$
\begin{bmatrix}A_0 & B_V \\ C_1 & D_V \end{bmatrix}
\begin{bmatrix} P_0^{-1} & 0 \\ 0 & I \end{bmatrix}
\begin{bmatrix}A_0^* & C_1^* \\ B_V^*  & D_V^* \end{bmatrix} =
\begin{bmatrix} P_0^{-1} & 0 \\ 0 & I \end{bmatrix}.
$$
In particular, $P_0^{-1}-A_0P_0^{-1}A_0^*= B_VB_V^*$ and $A_0P_0^{-1}C_1^*+B_VD_V^*=0$.

Using \eqref{Lyap} and the two equations just derived, we have
\begin{align*}
&P_0^{-1}P_1-A_0P_0^{-1}P_1\alpha= P_0^{-1}P_1-A_0P_0^{-1}(A_0^*P_1+C_1^*\gamma)
\\
=&P_0^{-1}P_1-A_0P_0^{-1}A_0^*P_1-A_0P_0^{-1}C_1^*\gamma 
\\
=&P_0^{-1}P_1-(P_0^{-1}-B_VB_V^*)P_1-A_0P_0^{-1}C_1^*\gamma
\\
=&B_VB_V^*P_1-A_0P_0^{-1}C_1^*\gamma =B_VB_V^*P_1+B_VD_V^*\gamma=B_VB_W^*,
\end{align*}
since $B_W^*$ is defined by $B_W^*=B_V^*P_1+D_V^*\gamma$. \hfill$\Box$. 
\medskip

\setcounter{equation}{0}

\section{Appendix: some remarks on the Douglas-Shapiro-Shields factorization}

In this section we discuss further properties of the DSS factorization, which first appeared in Section \ref{sec:three}. We divide the section into three parts, Part 6.1, Part 6.2 and Part 6.3.

\medskip

{\bf Part 6.1.} We start with $\Xi(z)=V(z)W^*(z)$, where $V(z)$ and $W(z)$ are rational bi-inner $m\times m$ matrix functions.  Furthermore, as in \eqref{realiV} and \eqref{realiW}, we  assume that
\[
V(z)=D_V+zC_V(I_{n_-}-zA_V)^{-1}B_V, \ W(z)=D_W+zC_W(I_{n_+}-zA_W)^{-1}B_W
\]
are bi-inner and stable. Moreover, both realizations are unitary realizations. Thus
$$
\begin{bmatrix} A_V & B_V \\ C_V & D_V\end{bmatrix} \mbox{\ and\ }
\begin{bmatrix} A_W & B_W \\ C_W & D_W\end{bmatrix}
$$
are both unitary matrices. These realizations are then minimal. Also, let $X$ be the unique solution 
of the Stein equation 
\begin{equation}\label{Steina}
X=A_VXA_W^*+B_VB_W^*.
\end{equation}
Then, see \cite{FK},
\begin{equation}\label{realizationXialt}
\Xi(z)=zC_V(I_{n_-}-zA_V)^{-1}B_R +{\Xi}_0 + C_R(zI_{n_+}-A_W^*)^{-1}C_W^*,
\end{equation}
where $B_R$, $C_R$  and  ${\Xi}_0$ are given by
\begin{align*}
B_R& = B_VD_W^*+A_VXC_W^* ,\\
C_R&=D_VB_W^* +C_VXA_W^*,\\
{\Xi}_0&=D_VD_W^*+C_VXC_W^*.
\end{align*}

Minimality of the realizations of $V(z)$ and $W(z)$ implies that the pairs $(C_V,A_V)$ and $(C_W,A_W)$ are observable.
Hence to show that the realization of ${\Xi}(z)$ given by \eqref{realizationXialt} is minimal, it suffices to prove the following proposition.

\begin{proposition}
The pair $(A_V,B_R)$ is controllable, and the pair $(C_R, A_W^*)$ is observable.
\end{proposition}

\noindent
{\bf Proof.} Assume the pair $(A_V,B_R)$ is not controllable. By the Hautus criterion there is an eigenvalue
$\lambda$ of $A_V^*$ and a non-zero vector $x$ such that $A_V^*x=\lambda x$ and
$B_R^*x=0$. Now $B_R^*=C_WX^*A_V^*+D_WB_V^*$. Thus
\begin{equation}\label{eq:basiscid}
C_WX^*A_V^*x+D_WB_V^*x=0.
\end{equation}
We shall use that 
$\begin{bmatrix} A_W & B_W \\ C_W & D_W \end{bmatrix}$ is unitary,
i.e.,
\begin{equation}\label{eq:Wunitary}
\begin{bmatrix}
A_W^* & C_W^* \\ B_W^* & D_W^*
\end{bmatrix}
\begin{bmatrix} A_W & B_W \\ C_W & D_W \end{bmatrix}
=\begin{bmatrix} I & 0 \\ 0 & I \end{bmatrix}.
\end{equation}
Firstly  multiply \eqref{eq:basiscid} on the left by $D_W^*$, to obtain that
\begin{equation*}
0  = D_W^*C_WX^*A_V^*x+ D_W^*D_WB_V^*x.
\end{equation*}
Using \eqref{eq:Wunitary} and \eqref{Steina} this gives
\begin{align*}
0&=-B_W^*A_WX^*A_V^*x+(I-B_W^*B_W)B_V^*x=
\\
&=-B_W^*(A_WX^*A_V^*+B_WB_V^*)x +B_V^*x=
\\
&= -B_W^*X^*x+B_V^*x.
\end{align*}
Secondly, multiply \eqref{eq:basiscid} by $C_W^*$, and use again \eqref{eq:Wunitary} and \eqref{Steina} to obtain that
\begin{align*}
0 & = C_W^*C_WX^*A_V^*x+ C_W^*D_WB_V^*x=
\\
&=(I-A_W^*A_W)X^*A_V^*x - A_W^*B_WB_V^*x=
\\
&=X^*A_V^*x-A_W^*(A_WX^*A_V^*+B_WB_V^*)x =
\\
&=X^*A_V^*x-A_W^*X^*x=(\lambda -A_W^*)X^*x.
\end{align*}
Combining with the fact that by assumption $A_V^*x=\lambda x$ we then obtain that
\begin{align*}
& \begin{bmatrix} B_V^* & -B_W^*\end{bmatrix}
\begin{bmatrix} x \\ X^*x \end{bmatrix} =0,
\\
& \begin{bmatrix} A_V^* & 0 \\ 0 & A_W^* \end{bmatrix}
\begin{bmatrix} x \\ X^*x \end{bmatrix} =\lambda \begin{bmatrix} x \\ X^*x \end{bmatrix}.
\end{align*}
However, the pair $\left( \begin{bmatrix} A_V & 0 \\ 0 & A_W \end{bmatrix},
\begin{bmatrix} B_V \\ -B_W\end{bmatrix}\right)$ is controllable since the realizations for $V(z)$ and $W(z)$ are minimal. So by duality and the Hautus criterion we must have 
$\begin{bmatrix} x \\ X^*x \end{bmatrix} =0$, which contradicts the assumption that $x\not= 0$.
This proves the first part of the claim.

The second part can be proved directly analogously to the first part. But alternatively, a proof by duality is perhaps more instructive. In fact, consider ${\Xi}^*(z)=W(z)V^*(z)$, which is unitary as well as ${\Xi}(z)$. The corresponding Stein equation is given by
$$
Y=A_WYA_V^* +B_WB_V^*.
$$
The unique solution is then given by $Y=X^*$. Analogously to the proof of the controllability of $(A_V, B_R)$ one show that the pair 
$(A_W, B_{R^*})$ is controllable. Observe that
$$
B_{R^*}=B_WD_V^*+A_WX^*C_V^*=C_R^*.
$$
Hence the pair $(A_W,C_R^*)$ is controllable, and thus  the dual pair $(C_R,A_W^*)$ is observable.
\hfill $\Box$

\medskip

{\bf Part 6.2.}  In this part we shall show that $V(z)$ and $W(z)$ are unique up to a unitary constant.

Suppose that ${\Xi}(z)=V(z)W^*(z)=V_1(z)W_1^*(z)$, where $V(z),W(z),V_1(z),\allowbreak W_1(z)$ take unitary values on the unit circle,
and are all rational bi-inner. Then, by Theorem 4.7.1 (iii) in \cite{FB}, we have that $V_1(z)=V(z)U$ and
$W_1(z)=W(z)U$ for some unitary constant $U$.

\medskip

Next, we consider how the solution of \eqref{Steina} changes when we consider the DSS factorization 
$\Xi(z)=V_1(z)W_1^*(z)$ instead of $\Xi(z)=V(z)W^*(z)$. 
In addition to \eqref{realizationXialt} we have
\begin{equation}\label{realizationXialt2}
{\Xi}(z)=zC_{V_1}(I_{n_-}-zA_{V_1})^{-1}\widehat{B_R} +R_0 + \widehat{C_R}(zI_{n_+}-A_{W_1}^*)^{-1}C_{W_1}^*,
\end{equation}
where
\begin{align*}
\widehat{B_R}& = B_{V_1}D_{W_1}^*+A_{V_1}X_1C_{W_1}^* ,\\
\widehat{C_R}&=D_{V_1}B_{W_1}^* +C_{V_1}X_1A_{W_1}^*,\\
R_0&=D_{V_1}D_{W_1}^*+C_{V_1}X_1C_{W_1}^*,
\end{align*}
and $X_1$ is the solution to the equation
$$
X_1=A_{V_1}X_1A_{W_1}^*+B_{V_1}B_{W_1} ^*.
$$
Because of the minimality of both realizations \eqref{realizationXialt} and \eqref{realizationXialt2}
there is an invertible matrix $S$ such that
\begin{equation}\label{eq:avav1}
A_{V_1}=S^{-1}A_VS, \ \ C_{V_1}=C_VS, \ \ B_{R}=S^{-1}\widehat{B_R}
\end{equation}
and there is an invertible matrix $T$ such that 
\begin{equation}\label{eq:awaw1}
A_{W_1}^*=T^{-1}A_W^*T, \ \ C_{W_1}^*=T^{-1}C_W^*, \ \  \widehat{C_R}=C_RT.
\end{equation}
Note that $S$ and $T$ are the state space similarities between the realizations $(A_V, B_R, C_V)$ and $(A_{V_1},\widehat{B_R}, C_{V_1})$, and $(A_W^*,C_R,C_W^*)$ and $(A_{W_1},\widehat{C_R},C_{W_1}^*)$, respectively.

Introduce $\Gamma_V={\rm col\,}(C_VA_V^j)_{j=0}^\infty$ and similarly 
$\Gamma_W={\rm col\,}(C_WA_W^j)_{j=0}^\infty$, $\Gamma_{V_1}={\rm col\,}(C_{V_1}A_{V_1}^j)_{j=0}^\infty$ and $\Gamma_{W_1}={\rm col\,}(C_{W_1} A_{W_1}^j)_{j=0}^\infty$.  Denote the block Toeplitz operator with symbol $\Xi(z)$ by $T_\Xi$.

Using formulas \eqref{eq:avav1} and \eqref{eq:awaw1}, and the formula $X=\Gamma_V^*T_\Xi\Gamma_W$ (see \cite {FK}),
 the following proposition is immediate.

\begin{proposition}
We have  $X=S^{-*}\Gamma_{V_1}^*T_\Xi\Gamma_{W_1}T^*$, 
and hence the matrices $X$ and $X_1=\Gamma_{V_1}^*T_\Xi\Gamma_{W_1}$ are related by $X_1=S^*XT^{-*}$. 
 \end{proposition}

\medskip

{\bf Part 6.3.} In this part we discuss the uniqueness of the unitary realizations of $V(z)$ and $W(z)$.

In comparison with Part 6.2 we further restrict the realizations of $V(z)$, $V_1(z)$, $W(z)$ and $W_1(z)$ to unitary realizations, that is,
$$
\begin{bmatrix} A_V & B_V \\ C_V & D_V\end{bmatrix} \mbox{\ and\ }
\begin{bmatrix} A_W & B_W \\ C_W & D_W\end{bmatrix}
$$
as well as
$$
\begin{bmatrix} A_{V_1} & B_{V_1} \\ C_{V_1} & D_{V_1}\end{bmatrix} \mbox{\ and\ }
\begin{bmatrix} A_{W_1} & B_{W_1} \\ C_{W_1} & D_{W_1}\end{bmatrix}
$$
are unitary, then the invertible matrices $S$ and $T$ in Part 6.3 are further restricted to being unitary. Indeed, 
$$
\begin{bmatrix}  A_{V_1} & B_{V_1} \\ C_{V_1} & D_{V_1}\end{bmatrix}=
\begin{bmatrix} S^{-1} & 0 \\ 0 & I \end{bmatrix}\begin{bmatrix} A_V & B_V \\ C_V & D_V\end{bmatrix}\begin{bmatrix} S & 0 \\ 0 & I \end{bmatrix}.
$$
Then, as $\begin{bmatrix}  A_{V_1} & B_{V_1} \\ C_{V_1} & D_{V_1}\end{bmatrix}$ is unitary, we obtain
$$
\begin{bmatrix}  A_{V}^* & C_{V}^* \\ B_{V}^* & D_{V}^*\end{bmatrix}
\begin{bmatrix} (SS^*)^{-1} & 0 \\ 0 & I \end{bmatrix}
\begin{bmatrix}  A_{V} & B_{V} \\ C_{V} & D_{V}\end{bmatrix} =
\begin{bmatrix} (SS^*)^{-1} & 0 \\ 0 & I \end{bmatrix}.
$$
Consider the $(1,1)$-block entry in this equation:
$$
A_V^*(SS^*)^{-1}A_V+C_V^*C_V=(SS^*)^{-1}.
$$
This means that $(SS^*)^{-1}$ is a solution to the matrix equation 
$$
P-A_V^*PA_V=C_V^*C_V.
$$
However, this equation has a unique solution, since $A_V$ is stable, and since $\begin{bmatrix} A_V & B_V \\ C_V & D_V\end{bmatrix} $ is unitary, the solution is equal to the identity. Hence $SS^*=I$, and so $S$ is unitary. Likewise, also $T$ is unitary.

It follows that the matrix $X$, when restricting the attention to unitary realizations of $V(z)$ and $W(z)$, is unique up to unitary equivalence.

This completes Part 6.3 and hence Section 6 is completed too. \hfill $\Box$

\medskip

\noindent
{\bf Acknowledgements.} The authors are very grateful to Jaap Korevaar for private communications which played a role in the earlier 2017 paper of the authors in \cite {GKR}.

The authors would also like to thank the referee for his/her remarks, which led to a substantial improvement of the paper. 

This work of the third author is based on research supported in part by the National Research Foundation of South Africa (Grant Number 145688).


\begin{thebibliography}{ww}
\bibitem{BGR}
J.A. Ball, I. Gohberg, L. Rodman, \emph{Interpolation of rational matrix functions.} Oper. Theory Adv. Appl.  45,  Birkh\"auser Verlag, Basel, 1990.

\bibitem{Ball-Ran1}
J.A. Ball, A.C.M. Ran, Local inverse spectral problems for rational matrix functions. \emph{Integr. Equ. Oper. Theory} 10, 349--415, 1987.


\bibitem{Ball-Ran2}
J.A. Ball, A.C.M. Ran, Global inverse spectral problems for rational matrix functions. \emph{Linear Algebra Appl.} 86, 237--282, 1987.


\bibitem{BGKOT21}
H. Bart, I. Gohberg, M.A. Kaashoek,  Invariants for Wiener-Hopf equivalence of analytic operator functions. In: \emph{Constructive methods of Wiener-Hopf factorization.} Oper. Theory Adv. Appl. 21, Birkh\"auser Verlag, Basel, 1986, 317-355.

\bibitem{BGKROT178}
H. Bart, I. Gohberg, M.A. Kaashoek, A.C.M. Ran,
\emph{Factorization of matrix and operator functions: the state space method}. 
Oper. Theory Adv. Appl.  178, Birkh\"auser, Basel, 2008.

\bibitem{BGKROT200}
H. Bart, I. Gohberg, M.A. Kaashoek, A.C.M. Ran, \emph{A state space approach to canonical factorization with applications}. Oper. Theory Adv. Appl.  200, Birkh\"auser, Basel, 2010.


\bibitem{CGOT3}
K. Clancey, I. Gohberg, \emph{Factorization of matrix functions and singular integral operators.}
Oper. Theory Adv. Appl.  3, Birkh\"auser Verlag, Basel, 1981.


\bibitem{FB}
A.E. Frazho, W. Boshri, \emph{An operator perspective on signals and systems.}  Oper. Theory Adv. Appl. 204, Birkh\"auser Verlag, Basel, 2010.

\bibitem{FK}
A.E. Frazho, M.A. Kaashoek, Canonical factorization of rational matrix functions. A note on a paper by P. Dewilde, \emph{Indag. Math.} 23 (2012), 1154--1164.

\bibitem{FKR}
A.E. Frazho, M.A. Kaashoek, A.C.M. Ran, The non-symmetric discrete algebraic Riccati equation and canonical factorization of rational matrix functions on the unit circle, \emph{Integr. Equ. Oper. Theory} 66 (2010), 215--229.

\bibitem{FR}
A.E. Frazho,  A.C.M. Ran, A note on inner-outer factorization of wide matrix-valued functions . In: \emph{Operator Theory, Analysis
and the State Space Approach}  Oper. Theory Adv. Appl. 271,  Birkh\"auser Verlag, Basel, 2018, 201--214.

\bibitem{GF}
I.C. Gohberg, I. A. Feldman, \emph{Convolution equations and projection methods for their solution.}
Transl. Math. Monogr., Vol. 41, Amer. Math. Soc., Providence, Rhode Island, 1974.

\bibitem{GGKOT49}
I. Gohberg, S. Goldberg, M.A. Kaashoek, \emph{Classes of linear operators}, Vol. I. OT 49,
Birkh\"auser Verlag, Basel, 1990.


\bibitem{GGKOT63}
I. Gohberg, S. Goldberg, M.A. Kaashoek, \emph{Classes of linear operators}, Vol. II. OT 63,
Birkh\"auser Verlag, Basel, 1993.



\bibitem{GKR}
G. Groenewald, M.A. Kaashoek, A.C.M. Ran, Wiener-Hopf indices of unitary functions on the unit circle in terms of realizations and related results on Toeplitz operators. \emph{Indag. Math.} 28 (2017) 694--710.

\bibitem{Potapov-60} V.P. Potapov, The multiplicative structure of $J$-contractive matrix functions. Tr. Mosk. Mat. Obs. 4 (1955), 125-236 (Russian); English translation: Amer. Math. Soc. Transl. (2) 15 (1960),  131-243.

\end{thebibliography}
\end{document}